\documentclass[12pt]{article}
\usepackage{amssymb,latexsym, amsmath, amscd, array, graphicx}
\usepackage{amssymb}

\usepackage{setspace}
\input xy
\xyoption{all}
\makeatletter 

\usepackage{setspace}
\usepackage{enumitem}

\date{}

\newtheorem{theorem}{Theorem}[section]

\newtheorem{proposition}[theorem]{Proposition}

\newtheorem{problem}[theorem]{Problem}

\newcommand{\z}{{\Bbb Z}}
\newcommand{\q}{{\Bbb Q}}

\newcommand{\N}{{\Bbb N}}

\newcommand{\lo}{\rightarrow}

\newcommand{\p}{\dim_{\z[\frac{1}{p}]} }

\newcommand{\black}{{\blacksquare}}

\begin{document}
\title{On actions of abelian Cantor groups}

\author{  Michael Levin\footnote{This research was supported by 
THE ISRAEL SCIENCE FOUNDATION (grant No. 522/14) }}

\maketitle
\begin{abstract} 
By  a Cantor group we mean a topological group homeomorphic 
to the Cantor set. The author earlier proved  
   that every compact metric space of
rational cohomological dimension $n$ can be obtained as 
the orbit space of a Cantor group action on a metric 
compact space of covering dimension $n$. In this paper, we consider actions of abelian
Cantor groups and, in particular,  
we show that in the  result  mentioned above
 the Cantor group can be assumed to be abelian for $n \geq 2$.
\bigskip
\\
{\bf Math. Subj. Class.:}  55M10, 22C05 (54F45)
\end{abstract}
\begin{section}{Introduction}\label{intro}
Throughout this paper we assume that spaces are separable metrizable,  actions
of topological groups are continuouus, and maps between spaces are continuous.
 We recall that a  compactum
means a compact metric space.

Let  $G$ be  an abelian group. The cohomological dimension $\dim_G X $ of a space $X$
is the smallest integer $n$ such that  the Cech cohomology $H^{n+1} (X, A; G)$ 
vanishes for every closed subset $A$ of $X$ if such an $n$ exists and $\dim_G X =\infty $ otherwise.
Clearly $\dim_G X \leq \dim X$ for every abelian group $G$ where $\dim X$
stands for the covering dimension of $X$.

 Exploring connections
between cohomological and covering dimensions is one of the central topics
in Dimension Theory. In this paper we focus on   connections related to actions of abelian Cantor groups.
By a Cantor group we mean a topological group homeomorphic to the Cantor set.  The following theorem was proved 
in \cite{levin-cantor}.

\begin{theorem}{\rm (\cite{levin-cantor})}
\label{levin-cantor}
Let $X$ be a compactum with $\dim_\q X =n$. Then there are an $n$-dimensional
 compactum  $Z$
and an action of a Cantor group  $\Gamma$ on $Z$ such that $X=Z/\Gamma$.
Moreover, the action of $\Gamma$ can be assumed to be free if $n=1$.
\end{theorem}

 It is known that for an action of a Cantor group $\Gamma$ on a compactum $Z$ we have 
 $\dim_\q Z=\dim_\q Z/\Gamma$.
  Thus  Theorem \ref{levin-cantor} provides a characterization of rational cohomological dimension in
  terms of Cantor group actions. The major  open problem regarding Theorem \ref{levin-cantor} is
  
  \begin{problem}
  \label{major-problem}
 Can  the Cantor group $\Gamma$ in Theorem \ref{levin-cantor}   be assumed to be  abelian?
 \end{problem}

  Actions of abelian Cantor groups were studied by Dranishnikov and  West in
  \cite{dranish-west}. 
  For a prime $p$ we denote by $\z^\N_p$ the Cantor group which is 
 the product of countably many copies of the $p$-cyclic
group $\z_p=\z/p\z$.  Dranishnikov and West  obtained the following important

\begin{theorem}
{\rm (\cite{dranish-west})}
\label{theorem-dranishnikov-west}
Let $X$ be a compactum. Then for every prime $p$
there are a compactum $Z$ and an action of the group 
$\Gamma=\z_p^\N$ on $Z$ such that
$\dim_{\z_p} Z\leq 1$ and $X =Z/\Gamma$.
\end{theorem}

Let $\cal P$ stand for the set of all prime numbers and let   $\z_* $  stand for
the Cantor group which is   the product
$\z_*=\Pi_{p \in \cal P} \z_p$. Denote by  $\z_*^\N$ the product of countably many copies
of $\z_*$.
It was pointed out in \cite{levin-cantor} that Theorem \ref{theorem-dranishnikov-west} implies

	\begin{theorem}{\rm (\cite{levin-cantor,dranish-west})}
	\label{derived-dranishnikov-west}
	Let $X$ be a compactum with $\dim_\q X \leq  n, n\geq 2$. Then there are
	a compactum $Z$  and an action of the  group $\Gamma=\z_*^\N$
	on $Z$ such that $\dim_\z Z \leq n $ and $X =Z/\Gamma$.
	\end{theorem}

We will refine both Theorem \ref{levin-cantor} and Theorem \ref{derived-dranishnikov-west}
by showing that
\begin{theorem}
\label{abelian-cantor}
Let $X$ be a compactum with $\dim_\q X \leq  n, n\geq 2$. Then there are
	a compactum $Z$  and an action of the  group $\Gamma=\z_*^\N$
	on $Z$ such that $\dim  Z \leq n $ and $X =Z/\Gamma$.
	\end{theorem}
	
	In order to prove Theorem \ref{abelian-cantor} we also need to refine 
	Theorem \ref{theorem-dranishnikov-west}.  Recall that a CW-complex is said to
	be an absolute extensor for a space $X$ if every map $f : A \lo K$ from
	a closed subset $A$ of $X$  continuously extends over $X$.
It is well-known that the cohomological dimension $\dim_G$ admits the following 
characterization: $\dim_G X \leq n$ if and only the Eilenber-MacLane complex $K(G,n)$
is an absolute extensor for $X$. Another important and well-known connection, known as Dranishnikov's
extension criterion,  between absolute extensors and
cohomological dimension was established in \cite{dranish-extension}.

\begin{theorem}{\rm (\cite{dranish-extension})}
\label{dranishnikov-extension}
${}$ 

(i) Let a CW-complex $K$ be an absolute extensor for a compactum $X$. Then
 $\dim_{H_n(K)} X \leq n$ for every $n$ where $H_n(K)$ is the reduced homology of $K$.
 
 (ii) Let  $X$ be  a finite dimensional compactum and $K$ 
a simply connected CW-complex $K$ such that 
$\dim_{H_n(K)} X \leq n$ for every $n$. Then $K$ is an absolute extensor for $X$.
\end{theorem}

Theorem \ref{dranishnikov-extension} was extended in \cite{dydak-levin-moore} to

\begin{theorem} {\rm (\cite{dydak-levin-moore})}
\label{dydak-levin-moore}
The standard Moore space $M(\z_p, 1)$ is an absolute extensor for finite dimensional compacta $X$ 
with $\dim_{\z_p} X \leq 1$.
\end{theorem}

The refinement of Theorem \ref{theorem-dranishnikov-west} that we need is 
\begin{theorem}
\label{refinement-dranishnikov-west}
Let $X$ be a compactum. Then for every prime $p$
there are a compactum $Z$ and an action of the group 
$\Gamma=\z_p^\N$ on $Z$ such that 
 $X =Z/\Gamma$ and $M(\z_p, 1)$  is an absolute extensor for $Z$.
\end{theorem}

Note that Theorem \ref{abelian-cantor} does not cover the case $n=1$  of
Theorem \ref{levin-cantor}. For this case we get only partial results.

\begin{theorem}
\label{abelian-dim=n}
Let $X$ be a compactum with $\dim_\q  X =1$ and $\dim  X=n$. Then there are
a compactum $Z$ and a free  action of an abelian Cantor group $\Gamma$ on $Z$  such that
$X=Z/\Gamma$ and 

(i) $\dim Z \leq n-1$ if $2\leq n\leq 6$ and 

(ii) $\dim Z \leq \big[ n/2\big]+2 $ if $n\geq 7$.
\end{theorem}

The paper is organized as follows: the proofs   are presented in
the following section and  few  comments are given in the last section.

\end{section}

\begin{section}{ Proofs}
\begin{proposition}
\label{epsilon-map}
Let $K$ be a CW-complex and $X$ a compactum such that for every $\epsilon>0$ there is
an $\epsilon$-map $g : X \lo Y$ to a compactum $Y$ such that $K$ is an absolute extensor for $Y$.
Then $K$ is an absolute extensor for $X$.
\end{proposition}
{\bf Proof.}
Let $f : A_X \lo K$ be a map from a closed subset $A_X$ of $X$. Consider a finite subcomplex
$L$ of $K$ such that $f(A_X) \subset L$. Embed $L$ into a Euclidean space $E$ and let 
$r : U \lo L$ be a retraction from a neighborhood $U$ of $L$ in $E$. Define $\delta>0$
to be such that for every $l \in L$ the closed $\delta$-ball $B(l, \delta)$ centered at $l$
is contained in $U$ and define $\epsilon>0$ to be such that for every subset  of $A_X$ of 
diameter$< \epsilon$ we have that the image of this subset under $f$ is of diameter$<4 \delta$
with respect to the Euclidean metric from $E$.

Let $g : X \lo Y$ be an $\epsilon$-map to a compactum $Y$ such that $K$ is an absolute extensor for $Y$.
Denote $A_Y=g(A_X)$ and take an open cover ${\cal V} $ of $Y$ such that
the cover $g^{-1}({\cal V})=\{ g^{-1}(V): V \in {\cal V} \}$ consists of
 sets of diameter$<2\epsilon$. Consider a partition of unity $\{ \phi_V | \phi _V: Y \lo [0,1], V \in {\cal V} \}$
 subordinate to ${\cal V}$,  fix a point $x_V \in g^{-1}(V)\cap A_X$
 for every  $V \in {\cal V}$ such that $V \cap A_Y\neq \emptyset$
 and 
 define the map $\phi : A_Y  \lo U$ by $\phi (y)= \sum_{V \in {\cal V}} f(x_V)\phi_V(y), y \in A_Y$.
 
 Note that $f$ and the map $\phi \circ g $ restricted to $A_X$  are homotopic as maps to $U$ and hence 
 $f$ and  the map $r\circ \phi \circ g$ restricted to $A_X$
   are homotopic as maps to $L$. Recall that $K$ is an absolute extensor for
 $Y$ and extend $r\circ \phi : A_Y \lo L$ to a map
 $f_Y : Y \lo K$. Then  $f$ and  $f_X = f_Y \circ g$ are homotopic on $A_X$ and hence $f$ extends over $X$.
 $\black$
 \\
 \\
 For proving  Theorem \ref{refinement-dranishnikov-west} we need the following result 
 by Dranishnikov and Uspenskij.
 \begin{theorem}
{\rm (\cite{d-u})}
\label{d-u}	
	Let $f : Z  \lo X$ be a 0-dimensional map  of compacta $Z$ and $X$. If   a CW-complex $K$
	 is an absolute extensor for $X$ then $K$ is an absolute extensor for $Z$ as well. In particular,
	 $\dim_G Z \leq \dim_G X $ for every abelian group $G$. 
 
\end{theorem}	
 {\bf Proof of Theorem \ref{refinement-dranishnikov-west}.}
 Let $\epsilon_n=1/n$ and  $ f_n : X \lo Y_n$  an $\epsilon_n$-map to a finite dimensional compactum $Y_n$.
 By Theorem \ref{theorem-dranishnikov-west} there is a compactum $Z^Y_n$ and
 an action of $ \z^{\N}_p$ on $Z^Y_n$ such that $\dim_{\z_p} Z^Y_n \leq 1$ and
 $Y_n=Z^Y_n / \z^\N_p$. Since the projection $\pi^Y_n : Z^Y_n \lo Y_n$ is $0$-dimensional we have
 that $Z^Y_n$ is finite dimensional as well and hence, by Theorem \ref{dydak-levin-moore},
 the Moore space $M(\z_p, 1)$ is an absolute extensor for $Z^Y_n$.
 
 Let $Z^X_n $ be the pull-back space of the maps $f_n$ and $\pi^Y_n$, and 
 $\pi^X_n : Z^X_n \lo X$ and $f^Z_n : Z^X_n \lo Z^Y_n$ the pull-back maps
  of $\pi^Y_n$ and $f_n$ respectively. Endow $Z^X_n$ with
 the pull-back action of  $\z^\N_p$.  Define $Z$ to be the pull-back space of the maps
 $\pi_n^X : Z^X_n \lo X$  for all $n$, consider  $Z$ with the pull-back action 
 of the group  $\Gamma$ being the product countably many copies $(\z^\N_p)_n, n \in\N,$ of
 $\z^\N_p$ and
 let  $\pi_n : Z \lo Z_n^X$  and
 $\pi : Z \lo X= Z/\Gamma$ be  the projections.
 
 Let us  show that $Z$ satisfies the conclusions of the theorem. 
 Clearly $\Gamma$ is isomorphic to $\z^\N_p$. The only thing that we need to check
 is that $M(\z_p, 1)$ is an absolute extensor for $Z$.  Take $\epsilon>0$.
 Since $\pi: Z \lo X$ is $0$-dimensional there is $\delta>0$ such that for every
 closed  subset $F$ of $X$ such that $F$ is of diameter$< \delta$  we have
 that every component  of $\pi^{-1}(F)$ is of diameter$< \epsilon$. Consider $n$ such that
 the fibers of $f_n : X \lo Y_n$ are of diameter$< \delta$. Then for every  fiber $F$ 
 of 
 the map $f^Z_n \circ \pi_n :Z \lo Z^Y_n$ we have that $ \pi (F)\subset X $ is 
 of diameter$< \delta$ since otherwise
  $f_n(\pi(F))=\pi^Y_n(f^Z_n(\pi_n (F)))$ cannot be a singleton.
 Thus every component  of  $F$  must be of diameter$<\epsilon$. 
 
 Consider the monotone-light decomposition of the map
 $f^Z_n \circ \pi_n$ with $g  : Z \lo Y$  and $g_Y : Y \lo Y_n$  being
 monotone and light respectively. Then, by Theorem \ref{d-u},
 $M(\z_p, 1)$ is an absolute extensor for $Y$ and  $g$ is an $\epsilon$-map. Thus
 for every $\epsilon>0$ the compactum $Z$ admits an $\epsilon$-map to a compactum for which
 $M(\z_p, 1)$ is an absolute extensor and hence, by Proposition \ref{epsilon-map},
  the Moore space $M(\z_p, 1)$ is an absolute extensor for $Z$. $\black$
 \\
 \\
 For proving Theorem \ref{abelian-cantor} we need
  \begin{theorem}
 \label{dydak-union-theorem}
 {\rm (\cite{dydak-union})}
 Let $X$ be a separable metric space and let  $K$ and $L$ be  CW-complexes.
 If   $X=A \cup B$ is the union of  subspaces $A$ and $B$ such that
 $ K$ and $L$ are absolute extensors for $A$ and $B$ respectively then  the  join $K*L$
 is an absolute extensor for $X$. In particular, if  $K$ is an absolute extensor  for $X$ then the suspension
 $\Sigma K$ is an absolute extensor for $X$ as well.

 \end{theorem}

 \begin{proposition}{\rm (\cite{levin-cantor})}
\label{levin-rational-moore}
 Let  $X$ be a compactum with $\dim_\q X \leq n$.  Then $  M(\q, n)$ is an absolute extensor for $X$.
 \end{proposition}
 
\begin{proposition}{\rm (\cite{levin-p-adic-free})}
\label{moore-spaces}
A compactum  $X$  is finite dimensional if and only if there is natural number $n$ such that the Moore spaces 
$ M(\q, n)$ and
$ M(\z_p, n)$ for  all prime $p$ are absolute extensors for $X$.
\end{proposition}
{\bf Proof of  Theorem \ref{abelian-cantor}.}
By Theorem \ref{refinement-dranishnikov-west} for every prime $p$ there is a compactum $Y_p$ and an action
of $\z_p^\N$ on $Y_p$ such that $X=Y_p/\z^\N_p$ and $M(\z_p, 1)$ is an absolute extensor for $X$.
Let $\pi_p : Y_p \lo X$ be the projection and $Z$ the pull-back space of the maps $\pi_p$ for all $p$.
Then for  the induced action of $\Gamma=\z^\N_*=\Pi_{p \in {\cal P}} \z^N_p$ on $Z$ we have
that $X=Z/ \Gamma$. Clearly the projections  $\pi : Z \lo X$ and  $\pi^Z_p : Z \lo Y_p$ are 
$0$-dimensional and hence,
by Theorem \ref{d-u}, we have that $\dim_\q Z \leq n$ and $M(\z_p, 1)$ for all $p$ are 
absolute extensors for $Z$. Then, 
by Proposition \ref{levin-rational-moore},  the Moore space $M(\q, n)$ is an absolute extensor for $Z$ and,
by Proposition \ref{dydak-union-theorem}, the Moore spaces 
$M(\z_p, n)=\Sigma^{n-1} M(\z_p, 1)$ for all prime $p$   are absolute extensors for $Z$ as well.
Thus, by Proposition \ref{moore-spaces},  the compactum $Z$ is finite-dimensional. 
By (i) of Theorem \ref{dranishnikov-extension}, we have that $\dim_{\z_p} Z  \leq 1$ for every prime $p$.
Since   $\dim_\q Z \leq n$ and $n \geq 2$  the Bockstein theory implies that $\dim_\z Z \leq n$.
Recall that  $Z$ is finite dimensional and hence  $\dim Z =\dim_\z  Z\leq n$. $\black$
\\
\\
Now we  turn to Theorem \ref{abelian-dim=n}.
  We will need
 
 \begin{proposition}
 \label{rational-circle} {\rm (\cite{levin-cantor})}
 Let  $X$ be a compactum with  $\dim_\q X\leq 1$. Then for every map 
 $f : A \lo S^1$ from a closed subset $A$ of $X$ there exists $k$ such that 
 $f$ followed by a $k$-fold covering map $S^1\lo S^1$ of $S^1$  continuously extends
 over $X$.
 \end{proposition}
 
 \begin{theorem}
 \label{torunczyk}
 {\rm (\cite{torunczyk})}
 Let $f : Z \lo X$ be a $0$-dimensional map of finite dimensional compacta. Then
 there is a subset $B$ of $Z$ such that $\dim B \leq\big[\frac{1}{2}\dim X\big]$ and
  $f$ is $1$-to-$1$ on $Z \setminus B$.
 \end{theorem}
 
 \begin{proposition}
 \label{lifting}
 Let $X$ be a compactum and $f : A \lo S^1$ a map from a closed subset $A$ of $X $  such that
 $f$ followed by a $k$-fold covering map $g : S^1\lo S^1$  extends to a map $\phi : X\lo S^1$.
 Consider the pull-back space $Z$ of the maps $\phi$ and $g$ and let 
 $\psi : Z \lo S^1$ and $h : Z \lo X$  be the pull-back maps of $\phi$ and $g$ respectively.
 Then $h$ restricted to $h^{-1}(A)$ and followed by $f$ extends over $Z$.
 \end{proposition}
 {\bf Proof.}
 Consider a component $C_Z$ of $A_Z=h^{-1}(A)$ and let  
$C=h(C_Z)\subset A\subset X$.
Note that $\psi: Z \lo S^1$ restricted to $C_Z$  and followed by $g$ coincides
with $g \circ (f|C )\circ ( h | C_Z)$.  Then $\psi$ restricted to $C_Z$ and
the map  $f\circ (h|C_Z)$  are liftings   via $g: S^1 \lo S^1$ of the map $\phi \circ (h|C_Z)$ and
hence, since $C_Z$ is connected,  differ by a rotation of $S^1$ and, as a result, are homotopic.
Thus $\psi$ and $f \circ h$ are homotopic on every component of $A_Z$ and therefore are
homotopic on $A_Z$ as well and hence $f \circ (h|A_Z)$ extends over $Z$.$\black$
\\
\\
 We  first prove the second part of Theorem \ref{abelian-dim=n}.
 \\
 \\
{\bf Proof  of (ii) of Theorem \ref{abelian-dim=n}.} Consider a countable collection 
of partial maps $f_i : A_i \lo S^1$ from closed subsets $A_i$ of $X$ to the circle $S^1$ 
 such that for every
partial map $f : A \lo S^1$ from a closed subset $A$ of $X$ there is $f_i : A_i \lo S^1$
such that $A\subset A_i$ and $f_i$ restricted to $A$ is homotopic to $f$.
By Proposition \ref{rational-circle} for every $f_i : A_i \lo S^1$  there is $k_i$ such that
$f_i$ followed by the $k_i$-fold covering map $g_i : S^1 \lo S^1$  of  $S^1$ extends 
to $\phi_i : X \lo S^1$.

 Let $Z_i$ be the pull-back space of the maps
$\phi_i$ and $g_i$ and let  $\psi_i : Z_i \lo S^1$ and  $h_i : Z_i \lo X$ be
 the pull-back maps of $\phi_i$ and $g_i$ respectively.
 Consider $g_i$ as the quotient map to the orbit space 
of the free  action of $\z_{k_i}$ on  $S^1$ and endow $Z_i$ with the pull-back action of
$\z_{k_i}$.  Then $X=Z_i/\z_{k_i}$.  

Define $Z$ to be the pull-back space of the maps $h_i$ for all $i$ and let
$\pi_i: Z \lo Z_i$ and $\pi : Z \lo X$ be the projections. Endow $Z$ with the pull-back action
of $\Gamma=\Pi_{\i \in \N} \z_{k_i}$. Clearly the action of $\Gamma$ on $Z$ is free and
 $X =Z/\Gamma$.

Let us show that $\dim Z \leq [n/2]+2$. Aiming at a contradiction assume that
$\dim Z > [n/2]+2$. By Theorem \ref{torunczyk} there is a subset $B \subset Z$ such that
$\dim B \leq [n/2]$ and $\pi$ restricted to $Z\setminus B$ is $1$-to-$1$.
Replacing $B$ by a larger  $G_\delta$-set of the same dimension we may assume that
$Z \setminus B$ is $\sigma$-compact. Since $\dim Z > [n/2]+2$ we have that $\dim (Z \setminus B)\geq 2$
and hence we can find a closed subset $F_Z \subset Z$ with $\dim F_Z \geq 2$ on which
$\pi$ is $1$-to-$1$. 

Then there is a closed subset $A_Z \subset F_Z$ and a map $f_Z : A_Z \lo S^1$ such that
$f_Z $ does not extend over $F_Z$. Since $\pi$ is $1$-to-$1$ on $A_Z$, the  map $f_Z$ induces
the map $f : A=\pi(A_Z) \lo S^1$ such that $f_Z =f\circ (\pi|A_Z)$. Take a partial map
$f_i : A_i \lo S^1$ such that $A\subset A_i$ and $f_i$ restricted to $A$ is homotopic to $f$.
Then 
by Proposition \ref{lifting} the map 
$h_i$ restricted to $A^Z_i=\pi_i(A_Z)$ and followed by $f_i$ extends over $Z_i$ and
hence $h_i$ restricted to $A^Z_i$ and followed by $f$  extends over $Z_i$ 
to a map $\psi : Z_i \lo S^1$. Thus $\pi_i$  followed by $\psi$ provides an extension of 
$f_Z$ over $F_Z$ and hence $\dim F_Z \leq 1$. This contradiction proves the theorem. $\black$
\\
\\
The proof  (i) of Theorem \ref{abelian-dim=n} relies on slight adjustments of the  relevant 
 constructions and  results  from \cite{levin-p-adic-free}.

Let  $k $ be   a natural number.
Consider a $2$-simplex $\Delta$, denote by $\Omega(k)$ the mapping cylinder of a $k$-fold covering map
$\partial \Delta \lo S^1$ and refer to the domain   $\partial \Delta$ and the  range $S^1$ 
of this map  as the bottom and the top of  
$\Omega(k)$ respectively. 

Given an increasing  sequence 
  $ k_1, k_2,...$ of natural numbers  such that $k_{i+1}$ is divisible by $k_i$
we will construct
a compactum $Y$  to which we will refer  as the  rational surface determined 
determined by the sequence $k_i$. The compactum $Y$ is
  constructed as the inverse limit 
of $2$-dimensional finite simplicial complexes $\Omega_i$.
Set  $\Omega_0$ to be  a $2$-simplex $\Delta$. 
Assume that $\Omega_i$ is constructed. Take a sufficiently fine triangulation of 
$\Omega_i$ and in every  $2$-simplex $\Delta $ of $\Omega_i$  remove the interior of $\Delta$ and attach to
$\partial \Delta$ the mapping cylinder $\Omega({k_{i+1}})$
 by identifying the bottom of $\Omega({k_{i+1}})$ with
$\partial \Delta$.
  Define $\Omega_{i+1}$ to be  the space that obtained this way from $\Omega_i$   and 
define the bonding map
$\omega_{i+1} : \Omega_{i+1} \lo \Omega_n$ to be a map that sends  each mapping cylinder 
$\Omega({k_{i+1}})$
to the corresponding simplex   $\Delta$ that identifies the bottom of $\Omega({k_{i+1}})$ with 
$\partial \Delta$ and
sends the top of $\Omega({k_{i+1}})$ to the barycenter of $\Delta$.
 We additionally assume that in the construction of the rational surface $\Omega$
the triangulations of  $\Omega_i$ are so fine that   the images of the simplexes of 
$\Omega_j$  in 
$\Omega_i$, $i < j$ under the map 
$\omega^i_j=\omega_j \circ \dots \circ \omega_{i+1} : \Omega_j\lo \Omega_i$
 are of diameter$<1/2^{j-i}$. Clearly $\dim \Omega=2$.

\begin{proposition}
\label{maps-to-surfaces}
A compactum $X$ with $\dim X \leq n, n\geq 2$,  and $\dim_\q X=1$ admits an $(n-2)$-dimensional
map  into a rational surface.
\end{proposition}
The  proof of Proposition \ref{maps-to-surfaces} is similar to the one of  
Proposition 2. 1 of \cite{levin-p-adic-free} (we just need replace 
 Proposition 1.5 of  \cite{levin-p-adic-free}  by
Proposition \ref{rational-circle} of the present paper).

\begin{proposition}
\label{surface-resolvable}
Every  rational  surface can be obtained as the orbit space of a free action
of an abelian Cantor group on a $1$-dimensional compactum.
\end{proposition}
The proof of Proposition \ref{surface-resolvable} is similar to 
the proof of Proposition 2.5 of \cite{levin-p-adic-free}. 
\\
\\
{\bf Proof of (i) of Theorem \ref{abelian-dim=n}.}
 By Proposition \ref{maps-to-surfaces}
there is an $(n-2)$-dimensional map $f : X  \lo \Omega$ to a rational  surface $\Omega$. 
By Proposition \ref{surface-resolvable} there is a free  action of  an abelian 
 Cantor group $\Gamma$  on a $1$-dimensional compactum $Y$
such that $Y /\Gamma =\Omega$. Let $Z$ be the pull-back of $f$ and the projection of $Y$
to $\Omega$.  Then the projection of $Z$ to $Y$ is an $(n-2)$-dimensional  map and hence
 $\dim Z  \leq \dim Y  +(n-2)\leq n-1$ and 
for the pull-back action of $\Gamma$ on $Z $ we have that $Z/\Gamma=X$. 
$\black$.

\end{section}

\begin{section}{Remarks}
From the proof of the factorizations theorems in \cite{levin-unstable} one can derive the following
modification of Proposition \ref{moore-spaces}: if  a compactum $X$ is such  that $\p X \leq 1$ and
$M(\z_p, 1)$ is an absolute extensor for $X$  then $\dim X \leq 2$.  This property together  with
Theorem \ref{refinement-dranishnikov-west} imply
\begin{theorem}
\label{[1/p]}
Let $X$ be a compactum with $\p X =1$. Then there are a compactum $Z$ with $\p Z=1, \dim_{\z_p} Z =1$,
 $\dim Z \leq 2$  and an action of $\Gamma=\z^\N_p$ on $Z$ such that $X=Z/\Gamma$.
\end{theorem}
An example of a compactum $X$ satisfying the conclusions of Theorem \ref{[1/p]}  was constructed
 \cite{dranish-west}, Theorem \ref{[1/p]} shows that such examples are of generic nature.
 
 Let us end with a few  problems that seem to be difficult and  interesting. Most of these problems
 are special cases  of  Problem \ref{major-problem} for $n =1$.  
 \begin{itemize}

\item
 Does there exist a $3$-dimensional compactum that can be obtained as the orbit space of an action
 of an abelian Cantor group on a $1$-dimensional compactum?  
 \item
 Does there exist a $3$-dimensional compactum that can be obtained as the orbit space of a free action
 of an abelian Cantor group on a $1$-dimensional compactum?  
 
 \item
 Does there exist an infinite-dimensional  compactum that can be obtained as the orbit space of an action
 of an abelian Cantor group on a $1$-dimensional compactum?  
 
 \item Can Theorem \ref{[1/p]} be improved to provide $Z$ with $\dim Z =1$?
 
 \item Does there exist an action of a Cantor group   on a $2$-dimensional compactum  with
 the orbit space being hereditarily infinite dimensional?
 \end{itemize}
 We recall that a compactum is said to be hereditarily infinite dimensional  if  its closed subsets are either
 $0$-dimensional or infinite dimensional.

\end{section}

Michael Levin\\
Department of Mathematics\\
Ben Gurion University of the Negev\\
P.O.B. 653\\
Be'er Sheva 84105, ISRAEL  \\
 mlevine@math.bgu.ac.il\\\\
\end{document}